\documentclass[10pt,leqno]{article}

\usepackage{amsmath,amsfonts,amscd,amssymb,theorem}

\long\def\comment#1\endcomment{}

\makeatletter
\begingroup
\gdef\th@dotted{\normalfont\itshape
  \def\@begintheorem##1##2{%
        \item[\hskip\labelsep \theorem@headerfont ##1\ ##2.]}%
\def\@opargbegintheorem##1##2##3{%
   \item[\hskip\labelsep \theorem@headerfont ##1\ ##2\ (##3).]}}
\endgroup
\makeatother

\theoremstyle{dotted}

\newtheorem{theorem}{Theorem}[section]
\newtheorem{lemma}[theorem]{Lemma}

\newtheorem{prop}[theorem]{Proposition}

\makeatletter
\begingroup
\gdef\th@upshape{\normalfont
  \def\@begintheorem##1##2{%
        \item[\hskip\labelsep \theorem@headerfont ##1\ ##2.]}%
\def\@opargbegintheorem##1##2##3{%
   \item[\hskip\labelsep \theorem@headerfont ##1\ ##2\ (##3).]}}
\endgroup
\makeatother

\theoremstyle{upshape}

\newtheorem{defn}[theorem]{Definition}
\newtheorem{remark}[theorem]{Remark}

\makeatletter
\renewcommand{\subsection}{\@startsection{subsection}{2}{0pt}{-3ex
plus -1ex minus -0.2ex}{-2mm plus -0pt minus
-2pt}{\normalfont\bfseries}} \makeatother

\makeatletter
\@addtoreset{equation}{section}
\makeatother

\newcommand{\cntrct}                
{\hspace{2pt}\raisebox{1pt}{\text{$\lrcorner$}}\hspace{2pt}}

\newcommand{\proof}[1][Proof.]{\smallskip\noindent{\em #1}}
\def\endproof{\hfill\ensuremath{\square}\par\medskip}

\def\eqref#1{\thetag{\ref{#1}}}

\let\latexref=\ref
\def\ref#1{{\normalfont{\latexref{#1}}}}

\newcommand{\wt}{\widetilde}

\newcommand{\hdot}{{\:\raisebox{3pt}{\text{\circle*{1.5}}}}}
\newcommand{\calo}{{\cal O}}
\newcommand{\Z}{{\mathbb Z}}
\newcommand{\E}{{\cal E}}
\newcommand{\D}{{\cal D}}
\newcommand{\T}{{\cal T}}

\newcommand{\g}{{\mathfrak g}}
\newcommand{\MM}{{\bf\sf M}}
\newcommand{\LL}{{\bf\sf L}}
\renewcommand{\L}{{\cal L}}

\newcommand{\X}{{\mathfrak X}}
\newcommand{\Y}{{\mathcal Y}}
\newcommand{\gm}{{\mathbb{G}_m}}

\newcommand{\eps}{\varepsilon}

\newcommand{\gr}{\operatorname{\sf gr}}

\newcommand{\cchar}{\operatorname{\sf char}}
\newcommand{\Spec}{\operatorname{Spec}}
\newcommand{\spec}{\operatorname{{\cal S}{\it pec}}}

\newcommand{\Tot}{\operatorname{{\sf Tot}}}

\newcommand{\sllash}{/\!\!\,/}

\title{On the coordinate ring of a projective Poisson scheme}

\author{D. Kaledin\thanks{Partially supported by CRDF grant
RM1-2354-MO-02.}}

\date{}

\begin{document}

\maketitle


\section*{Introduction.}

Let $Y$ be a smooth symplectic algebraic variety over a field $k$ of
characteristic $0$, or, more generally, a Poisson scheme over
$k$. Assume that a reductive group $G$ acts on $Y$ preserving the
Poisson structure. Moreover, assume that the $G$-action is
Hamiltonian -- namely, that there exists a Lie algebra homomorphism
$\mu:\g \to H^0(Y,\calo_Y)$ from the Lie algebra $\g$ of the group
$G$ to the Lie algebra of global functions on $X$ such that any $\xi
\in \g$ acts on $\calo_Y$ by taking Poisson bracket with the
function $\mu(\xi)$. Then there exists a well-known procedure of
taking quotient with respect to the $G$-action, appropriate for the
category of Poisson schemes. The procedure is called {\em
Hamiltonian reduction}. It consists of two steps. Firstly, one
considers the map $\mu$ as a so-called {\em moment map} $\mu:Y \to
\g^*$ from $Y$ into the affine space $\g^*$ carrying the coadjoint
representation of $G$. Secondly, one replaces $Y$ with the zero
fiber $Y_0 = \mu^{-1}(0) \subset Y$ of the moment map and takes the
usual, suitably interpreted quotient $Y_0 \sllash G$ with respect to
the $G$-action (usually the quotient is taken in the sense of
GIT). If the group $G$ is abelian, one can also do the steps in the
opposite order: firstly, one takes the quotient $Y \sllash G$, then
takes the zero fiber $X = (Y \sllash G)_0$. In any case, the
resulting scheme $X = Y_0 \sllash G$ is automatically Poisson;
moreover, if $Y$ is symplectic and $X$ is smooth, then $X$ is
symplectic as well.

The goal of this note is to show how to invert the Hamiltonian
reduction procedure in the particular case when the group $G$ is the
multiplicative group $\gm$. We start with a Poisson scheme $X$ and a
line bundle $L$ on the scheme $X$. We construct, under certain
assumptions, a Poisson deformation $\X/S$ of the scheme $X$ over the
formal disc $S = \Spec k[[t]]$ and a line bundle $\L$ on $\X$
extending $L$. Then we prove that the total space $Y$ of the
corresponding $\gm$-principal bundle on $\X$ is a Poisson scheme,
and that the natural $G$-action on $Y$ is Hamiltonian, with the
projection $Y \to \X \to S$ giving the moment map. Thus $X$ is the
reduction of the Poisson scheme $Y$ with respect to a Hamiltonian
$\gm$-action. Moreover, we show that under additional assumptions
such a deformation $\X/S$ is essentially unique.

It seems natural to approach this problem by means of the general
Poisson deformation theory, as described in, for instance, the paper
\cite{GK}. However, it turns out that there exists a simpler
proof. In fact, our proof is almost elementary; the only slightly
non-trivial technique that we use is the notion of a stack of
groupoids. We generally follow the lines of the paper \cite{KV},
where a similar (but harder) result was proved for deformations of
smooth symplectic manifolds.

Thus, the question of inverting Hamiltonian reduction with respect
to $\gm$ admits an answer, which is at once simple and more or less
complete. Perhaps the main question the reader would ask at this
point is: why should he or she care? -- after all, the question
itself is dangerously close to pure sophistry. The author would like
to quote the usual meek excuse of ``the result being interesting in
its own right'', and add, as an explanation, that he arrived at the
problem while studying Poisson schemes that are
projective. Intuitively, the projective coordinate ring of such a
scheme ought to carry a Poisson structure, but it does not (and this
problem was noted already by M. Kontsevich in \cite{K}). Our theorem
explains why not: one has to add an additional deformation. This
might be useful in the study of quantizations of such schemes.

We were not able to extend our result to arbitrary reductive groups,
and it might be that such a generalization simply does not
exist. However, aside from certain cohomology vanishing, we impose
absolutely no restrictions on $X$: it may be singular, non-normal,
reducible, non-reduced.

\subsection*{Acknowledgments.} I am grateful to R. Bezrukavnikov,
D. Huybrechts and A. Kuznetsov for valuable discussions. I am
grateful to the referee for a detailed report and important
suggestions; in particular, he suggested that Remarks~\ref{prrd}
and~\ref{poi.md} would clarify the exposition.

\section{Statements and definitions.}

Fix a base field $k$ of characteristic $\cchar k = 0$. By a {\em
Poisson algebra} $A$ over $k$ we will understand a unital
commutative algebra $A$ over $k$ equipped with a skew-linear
operation $\{-,-\}:\Lambda^2_k(A) \to A$ such that
$$
\{a,bc\} = \{a,b\}c + \{a,c\}b, \quad
0 = \{a,\{b,c\}\} + \{b,\{c,a\}\} + \{c,\{a,b\}\},
$$
for all $a,b,c \in A$. The operation $\{-,-\}$ is called the {\em
Poisson bracket}. By the {\em Poisson center} of a Poisson algebra
$A$ we will understand the subspace of elements $a \in A$ such that
$\{a,A\} = 0$. By a {\em Poisson module} $M$ over the algebra $A$ we
will understand an $A$-module $M$ equipped with a Poisson bracket
operation $\{-,-\}:A \otimes_k M \to M$ such that
\begin{align*}
\{a,\{b,m\}\} &= \{b,\{a,m\}\} + \{\{a,b\},m\},\\
\{a,bm\} &= \{a,b\}m + b\{a,m\},\\
\{ab,m\}&=a\{b,m\}+b\{a,m\},
\end{align*}
for all $a,b \in A$, $m \in M$. By a {\em Poisson scheme} $X$ we
will understand a scheme $X$ equipped with a Poisson bracket in its
structure sheaf. By a {\em Poisson sheaf} over a Poisson scheme $X$
we will understand a coherent sheaf $\E$ of $\calo_X$-modules
equipped with a Poisson bracket. Note that the spectrum $\Spec A$ of
a Poisson algebra $A$ is an affine Poisson scheme, and the
localization $\wt{M}$ of a Poisson $A$-module $M$ is a Poisson sheaf
on $\Spec A$. By a {\em Poisson line bundle} we will understand a
Poisson sheaf which is flat (equivalently, locally trivial in
Zariski topology) and has rank $1$. Given a line bundle $\L$ on a
scheme $X$, we will denote by
$$
\Tot(\L) = \spec \bigoplus_{p \in \Z} \L^p
$$
the total space of $\L$ without a zero section. Then if $\L$ and $X$
are Poisson, the scheme $\Tot(\L)$ is also Poisson. If a Poisson
scheme $X$ is equipped with a map $f:X \to S$ to another scheme $S$,
we will say that $X$ is {\em Poisson over $S$} whenever the subsheaf
$f^{-1}\calo_s \subset \calo_X$ lies in the Poisson center of the
algebra sheaf $\calo_X$.

\begin{remark}\label{prrd}
Explictly, the Poisson bracket on $\Tot(\L)$ for a Poisson line
bundle $\L$ is defined as follows. Local functions on $\Tot(\L)$ are
linear combinations of local sections $s_n \in H^0(U,\L^{\otimes
n})$, $U \subset \X$, $n \in \Z$. It suffices to define
$\{s_n,s_m\}$, and it suffices to do it locally on $\X$. Choose a
trivialization $s \in H^0(U,\L)$ over some Zariski open $U \subset
\X$; then we have $s_n=fs^n$ and $s_m=gs^m$ for some functions $f,g
\in H^0(U,\calo_U)$. By definition, we must have
$$
\begin{aligned}
\{s_n,s_m\}&=\{fs^n,gs^m\}\\
&=\{f,g\}s^{n+m}+(mg\{f,s\}-nf\{g,s\})s^{n+m-1}+nmfg\{s,s\}s^{n+m-2}\\
&=\{f,g\}s^{n+m}+(mg\{f,s\}-nf\{g,s\})s^{n+m-1},
\end{aligned}
$$
where $\{f,s\},\{g,s\} \in H^0(U,\L)$ are defined by the Poisson
structure on the line bundle $\L$. It elementary to check that this
does not depend on the choice of the trivialization $s$.
\end{remark}

Let $X$ be a Poisson scheme. Any function $f$ on $X$ defines a
vector field $H_f$ on $X$ by setting
\begin{equation}\label{h_f}
H_f(g) = \{f,g\}
\end{equation}
for any (local) function $g \in \calo_X$. The vector field $H_f$ is
called the {\em Hamiltonian vector field} of the function $f$. It
vanishes if and only if the function $f$ lies in the Poisson center
of $\calo_X$. Moreover, for any Poisson sheaf $\E$ on $X$
\eqref{h_f} defines a differential operator $H_f:\E \to \E$ of order
$1$. If $f$ is central, so that $H_f = 0$ on $\calo_X$, then $H_f:\E
\to \E$ is a $\calo_X$-linear map.

We will say that a Poisson scheme $X$ is {\em non-degenerate} at a
smooth closed point $x \in X$ if the whole tangent space $T_x$ is
spanned by Hamiltonian vector fields $H_f$ corresponding to local
functions $f \in \calo_{X,x}$ on $X$ near $X$. If $X$ is smooth, we
will say that it is {\em symplectic} if it is non-degenerate at all
point. It is well-known that this is equivalent to saying that $X$
carries a non-degenerate closed $2$-form $\Omega \in
H^0(X,\Omega^2_X)$ (this is the standard definition nowadays, but
historically, the Poisson definition came first).

\begin{remark}\label{poi.md}
In the case when the Poisson scheme $X$ is non-degenerate, a Poisson
module is tautologically the same as a $\D$-module -- that is, a
sheaf of modules over the algebra $\D_X$ of differential operators
on $X$. In particular, the category of Poisson modules does not
depend on the symplectic structure. In general, the forgetful
functor from Poisson modules to coherent sheaves on $X$ admits a
left-adjoint, $\E \mapsto U_X \otimes \E$, where $U_X$ is a certain
algebra sheaf on $X$, a sort of an universal enveloping algebra for
the Poisson structure. Explicitly, the algebra $U_X$ is generated by
$\calo_X$ and $\Omega^1_X$, subject to relations $[a,(db)]=\{a,b\}$
and $[(da),(db)]=d\{a,b\}$. Just as the algebra $\D_X$, the algebra
$U_X$ has an increasing filtration by order; the associated graded
algebra $\gr U_X$ is isomorphic to the symmetric algebra
$S^\hdot(\Omega^1(X))$. If the Poisson structure is non-degenerate,
we have $\Omega^1(X) \cong \T(X)$ and $U_X \cong \D_X$. In the
general case the algebra $U_X$ is probably useless; if $X$ is
smooth, it might be convenient for some applications. However, we
will not need it, so we skip all the details.
\end{remark}

Denote by $S = \Spec k[[t]]$ the formal disc over $k$. For any
Poisson scheme $X$ equipped with a map $f:X \to S$, the pullback
$f^*t$ of the parameter $t$ gives a global function on $X$, and we
have the Hamiltonian vector field $H_t$. By our convention, $X$ is a
Poisson scheme over $S$ if and only if $H_t=0$. Another important
situation is when $H_t$ can be integrated to give a $\gm$-action on
$X$. In this case, we say that the $\gm$-action is Hamiltonian, and
that $f:X \to S$ is the {\em moment map} for the $\gm$-action. Being
a Hamiltonian vector field, $H_t$ is in particular a derivation with
respect to the Poisson bracket; therefore a Hamiltonian $\gm$-action
preserves the Poisson structure.

We leave it to the reader to check that all of the above make
sense {\em verbatim} if one allows $X$ to be a formal scheme.

\begin{defn}\label{ms}
A {\em formal moment system} is the pair $\langle \X/S, \L\rangle$
of a flat Poisson formal scheme $\X/S$ and a Poisson line bundle
$\L$ on $\X$ such that
\begin{enumerate}
\item the Hamiltonian vector field $H_t$ vanishes on $\calo_X$ ($X$
is Poisson over $S$), and
\item the $\calo_X$-linear map $H_t:\L\to \L$ is equal to the
identity map.
\end{enumerate}
\end{defn}

Note that since $\X$ is Poisson over $S$, the fiber $\X_o = \X
\times_S o$ of the scheme $\X$ over the closed point $o \in X$ is
also a Poisson scheme, and the line bundle $\L$ gives by restriction
a line bundle $L$ on $\X_o$. We will say that the moment system
$\langle \X/S,\L \rangle$ {\em extends} the pair $\langle \X_o,L
\rangle$. Our main result is the following. 

\begin{theorem}\label{main}
Assume given a Poisson scheme $X$ and a line bundle $\L$ on $X$.
\begin{enumerate}
\item Assume that $H^2(X,\calo_X) = 0$. Then there exists a formal
moment system $\langle \X/S,\L \rangle$ extending the pair $\langle
X, L \rangle$.
\item Assume in addition that $H^1(X,\calo_X) = 0$. Then such a
moment system is unique up to an isomorphism.
\end{enumerate}
\end{theorem}

In applications, one often has one additional piece of structure on
a Poisson scheme $X$ -- namely, a vector field $\xi$ which is {\em
conformal of weight $\lambda \in k$} in the sense that
\begin{equation}
\xi(\{f,g\}) = \{\xi(f),g\} + \{f,\xi(g)\} + \lambda\{f,g\}.
\end{equation}
Theorem~\ref{main} is compatible with such a vector field in the
following sense.

\begin{prop}\label{conf}
In the assumptions of Theorem~\ref{main}~\thetag{ii}, assume in
addition that $X$ is equipped with a conformal vector field $\xi$ of
some weight $\lambda$. Then $\xi$ extends to a vector field on $\X$
such that $\xi()=\lambda t$, and $\L$ admits an action of $\xi$
compatible with the Poisson module structure.
\end{prop}

To state our second result, say that a scheme $X$ is {\em
algebraically convex} (by analogy with holomorphically convex
spaces) if the natural morphism
$$
X \to \Spec H^0(X,\calo_X)
$$
is a projective map.

\begin{theorem}\label{alg}
In the assumptions of Theorem~\ref{main}, assume that in addition to
$H^1(X,\calo_X) = H^2(X,\calo_X) = 0$, the scheme $X$ is
algebraically convex. Then the formal scheme $\X$ can be extended to
an actual algebraically convex scheme which is flat over $S$. The
algebra $H^0(\X,\calo_\X)$ of global functions on $\X$ is flat over
$k[[t]]$ and $t$-adically complete, the natural map
$H^0(\X,\calo_X)/t \to H^0(X,\calo_X)$ is an isomorphism, the
natural $\gm$-action on the total space $\Tot(\L)$ is Hamiltonian,
and $\Tot(\L) \to \X \to S$ gives a moment map. If $X$ is a regular
scheme, then both $\X$ and $\Tot(\L)$ are regular schemes, while
$H^0(\X,\calo_\X)$ is a normal $k$-algebra. If $X$ is symplectic,
then $\Tot(\L)$ is symplectic.
\end{theorem}

\subsection{The symplectic case.}
Assume that the Poisson scheme $X$ is smooth and {\em symplectic},
so that the Poisson structure is defined by a symplectic form
$\Omega_X \in H^0(X,\Omega^2_X)$. For any deformation $\langle
\X,\Omega_{\X}\rangle$ of the pair $\langle X,\Omega_X \rangle$ over
a local Artin base $S = \Spec \calo_S$, the Gauss-Manin connection
provides a trivialization $H^\hdot_{DR}(\X/S) \cong H^\hdot_{DR}(X)
\otimes \calo_S$ of the relative de Rham cohomology module
$H^\hdot_{DR}(\X/S)$. The cohomology class $[\Omega_\X] \in
H^2_{DR}(\X/S)$ need not be constant with respect to this
trivialization; on the contrary, one defines the {\em period map}
$P:S \to H^2_{DR}(X)$ as the graph of $[\Omega_\X] \in
H^2_{DR}(\X/S) \cong H^2_{DR}(X) \otimes \calo_S$, and it has been
shown in \cite{KV} that the period map completely defines the
deformation (in particular, $P$ is constant if and only if $\X \cong
X \times S$ and $\Omega_{\X} = \Omega_X$). The moment systems are
described in this language in the following way.

\begin{lemma}\label{sympl}
Assume given a line bundle $L$ on a symplectic manifold $\langle
X,\Omega \rangle$, and a formal moment system $\langle
\X/S,\L\rangle$, $S = \Spec k[[t]]$. Then the Poisson scheme $\X/S$
is non-degenerate, and the corresponding relative symplectic form
$\Omega_{\X}$ on $\X/S$ satisfies
\begin{equation}\label{line}
[\Omega_{\X}] = [\Omega_X] + [L]t \in H^2_{DR}(\X/S)
\cong H^2_{DR}(X) \otimes \calo_S = H^2_{DR}(X)[[t]],
\end{equation}
where $[L] \in H^2_{DR}(X)$ is the first Chern class of the line
bundle $L$.
\end{lemma}

\proof{} By definition, for any smooth family $\X/S$ and a closed
relative differential form $\alpha \in H^0(\X,\Omega^\hdot_{\X/S})$,
one can compute the Gauss-Manin connection $\nabla[\alpha] \in
H^l_{DR}(\X/S) \otimes \Omega^1_S$ by the following procedure: one
lifts the relative form $\alpha$ to an absolute form $\wt{\alpha}
\in H^0(\X,\Omega^\hdot_{\X})$, one notices that $d\wt{\alpha} =
d\alpha = 0 \mod \Omega^{\geq 1}_S$, and, independently of the
choice of a lifting, one has $[d\wt{\alpha}]=\nabla[\alpha] \mod
\Omega^{\geq 2}_S$, where both sides are hypercohomology classes
with coefficients in
$$
\Omega^\hdot_{\X/S} \otimes \Omega^1_S \cong
\left(\Omega^\hdot_{\X}\cdot\Omega^{\geq 1}_S\right)/
\left(\Omega^\hdot_{\X}\cdot\Omega^{\geq 2}_S\right).
$$
Denote by $Y=\Tot(\L)$ the total space of the line bundle $\L$
without the zero section, and let $\pi:Y \to \X$ be the natural
projection. To compute $H^\hdot_{DR}(\X/S)$, one can use the
equivariant de Rham complex $\Omega^\hdot_{\gm}(Y)$ (see, e.g.,
\cite{BGV}). Recall that we have
$$
\Omega^\hdot_{\gm}(Y) = \Omega^\hdot(Y) \otimes k[u],
$$
where $u$ is an additional generator of degree $2$; the differential
$d^{\gm}$ in the equivariant de Rham complex is given by
$d^{\gm}(\alpha) = d\alpha + ui_\xi(\alpha)$, where $i_\xi$ is the
contraction with the differential $\xi \in H^0(Y,\T_Y)$ of the
$\gm$-action on $Y$. By definition, the natural map
$\pi^*\Omega^\hdot_{\X/S} \to \Omega^\hdot_Y$ extends to a map
$$
\pi^*\Omega^\hdot_{\X/S} \to \Omega^\hdot_{\gm}(Y),
$$
which induces an isomorphism on hypercohomology groups.

Now, the relative symplectic form $\pi^*\Omega_{\X} \in
H^0(Y,\Omega^2_{\gm}(Y))$ by definition comes from the absolute
symplectic form $\Omega_Y \in H^0(Y,\Omega^2_Y)$. The form
$\Omega_Y$ is not closed in the equivariant de Rham complex: while
$d\Omega_Y=0$, $i_\xi\Omega_Y$ is non-trivial. In fact, since the
projection $Y \to S$ is the moment map for the $\gm$-action on $Y$,
we have $i_\xi\Omega_Y = dt$. We conclude that
$$
\nabla[\Omega_{\X}]=\nabla[\pi^*\Omega_{\X}]=[d^{\gm}\Omega_Y]=[ui_\xi\Omega_Y]
= [u]dt.
$$
It remains to notice that by definition, we have $\pi^*[L] = u \in
H^2(Y,\Omega_{\gm}^\hdot(Y)) \cong H^2_{DR}(\X/S)$. Therefore
$\nabla[\Omega_\X]=[L]dt$, which, since $[\Omega_\X]=[\Omega_X] \mod
t$, implies \eqref{line}.
\endproof

To understand better the relation between symplectic deformations
and the moment systems, assume that $X$ is a symplectic manifold
with
$$
H^1(X,\calo_X) = H^2(X,\calo_X) = 0.
$$ 
Then the main theorem of \cite{KV} claims that $X$ admits a {\em
universal} symplectic formal deformation $\X/C$ over a smooth base
$C$. The corresponding period map $P:C \to H^2_{DR}(X)$ identifies
$C$ with the formal neighborhood of the cohomology class $[\Omega_X]
\in H^2_{DR}(X)$ in the affine space $H^2_{DR}(X)$. Then
Lemma~\ref{sympl} shows that a moment system associated with a line
bundle $L$ is parametrized by the line $[\Omega_X]+t[L] \in C$ in
the base $C$ of the universal deformation $\X/C$.

The methods of \cite{KV} work only in the symplectic situation;
moreover, they do not give a Poisson line bundle $\L$. On the other
hand, the approach of \cite{KV} requires weaker cohomological
assumptions on $X$ and, more importantly, it provides a
multi-parameter deformation $\X/C$ which is universal in appropriate
sense. However, for general Poisson schemes -- for instance, for a
Poisson scheme with zero Poisson bracket -- the deformation theory
behaves in an unpredictable way, and there is no reason to expect
that a universal Poisson deformation even exists, let alone is
smooth. So, in general, Theorem~\ref{main} is not far from an
optimal possible result on Poisson deformations.

\section{Proofs.}

We will now prove Theorem~\ref{main}, Proposition~\ref{conf} and
Theorem~\ref{alg}. Our proof proceeds by induction. To set up the
induction, for any integer $n \geq 1$ we denote by $S_n = \Spec
k[[t]]/t^{n+1}$ the $n$-th infinitesemal neighborhood of the special
point $o \in S$. By definition, we have canonical embeddings
$$
S_1 \subset \dots \subset S_n \subset \dots.
$$ 

\begin{defn}\label{ms.n} 
An {\em order-$n$ moment system} is a pair $\langle X_n/S_n, \L
\rangle$ of a flat Poisson scheme $X_n/S_n$ and a line bundle $\L$
on the subscheme $X_{n-1} = X_n \times_{S_n} S_{n-1} \subset X_n$
equipped with a structure of a Poisson sheaf on $X_n$, such that
$X_n/S_n$ and $\L$ satisfy the conditions \thetag{i}, \thetag{ii} of
Definition~\ref{ms}.
\end{defn}

Note that the condition \thetag{i} insures that $X_{n-1} \subset
X_n$ is a Poisson scheme, and $t^n\L \subset \L$ is a Poisson
subsheaf. Thus, given an order-$n$ moment system $\langle X_n/S_n,\L
\rangle$, we obtain by restriction a moment system $\langle X_{n-1},
\L/t^n\L\rangle$ of order $n-1$. Analogously, given a formal moment
system, we obtain by restriction an order-$n$ moment system for
every $n \geq 1$. For $n=0$, we simply obtain a Poisson scheme
$X_0$. We will say that the given moment system {\em extends} the
Poisson scheme $X_0$.

\begin{remark} 
The sheaf $\L$ in Definition~\ref{ms.n} is {\em not} a Poisson
module over $X_{n-1}$ (multiplication by $t^n$ is trivial on $\L$,
but the bracket with $t^n$ is not trivial -- in fact, up to a
constant $H_{t^n}$ is equal to the multiplication by $t^{n-1}$).
\end{remark}

Moment systems of order $n$ form a category, with morphisms (we will
only need isomorphisms) defined in the obvious way. For an arbitrary
Poisson scheme $X$, the pair $\langle X \times S_n, \calo_{X_{n-1}}
\rangle$ is obviously an order-$n$ moment system. We will call it
trivial. Given an order-$n$ moment system $\langle \X_n,\L \rangle$,
by a trivialization of the sheaf $\L$ we will understand a section
$e \in H^0(X_n,\L)$ which gives an isomorphism $\calo_{X_{n-1}} \to
\L$.

\begin{lemma}\label{triv}
Let $X = \Spec A$ be an affine Poisson scheme. Then for any
order-$n$ moment system $\langle X_n,\L\rangle$ extending $X$, every
trivialization $e$ of the sheaf $\L$ on $X_{n-1}$ extends uniquely
to an isomorphism between $\langle X_n,\L\rangle$ and the trivial
moment system $\langle X \times S_n, \calo_{X_{n-1}} \rangle$.
\end{lemma}

\proof{} Since $X_n$ extends an affine scheme $X$, it is itself
affine, so that $X_n = \Spec A_n$ for some Poisson algebra $A_n$. We
have $A_n/t \cong A$, and $A_{n-1}=A_n/t^n$ gives $X_{n-1} = \Spec
A_{n-1}$. The sheaf $\L$ corresponds to a $A_n$-module which is by
assumption identified with $eA_{n-1}$. Denote by
$$
A^0_n \subset A_n
$$
the subalgebra of functions $a \in A_n$ such that $\{a,e\}=0$. We
claim that the restriction map $A_n \to A_n/t \cong A$ induces an
isomorphism $A^0_n \cong A$. Indeed, this is obvious for the trivial
moment system. By induction on $n$, we may assume that it is true
for $A_{n-1}^0 \subset A_{n-1} = A_n/t^n$. Thus it suffices to prove
that for every function $a \in A_{n-1}^0$, the space $P_a \subset
A_n$ of all functions $b \in A$ such that $b = a \mod t^n$ contains
exactly one function $a' \in P_a$ satisfying $\{a',e\} = 0$. But
this is obvious: the space $P_a$ is a torsor over $t^nA \subset
A_n$, the commutator $\{P_a,e\} \subset eA_{n-1}$ lies in
$t^{n-1}eA_{n-1} \subset eA_{n-1}$, and the commutator map
$\{-,e\}:A_n \to eA_{n-1}$ identifies $t^nA \subset A_n$ with
$t^{n-1}eA_{n-1}$.

To finish the proof, it suffices to notice that $A_0 \otimes
k[t]/t^{n+1} \to A_n$ is a Poisson map, hence an isomorphism between
$A_n$ and the trivial moment system $A \otimes k[t]/t^{n+1}$.
\endproof

Let now $X$ be an arbitrary Poisson scheme, and assume given a
moment system $\langle X_n,\L \rangle$ of order $n$ which extends
$X$. Denote by $\LL_n$ the groupoid of all line bundles $\L'$ on
$X_n$ equipped with an isomorphism $\L'/t^n \cong \L$. Denote by
$\MM_n$ the groupoid of all moment systems $\langle X_{n+1},\L
\rangle$ of order $(n+1)$ equipped with an isomorphism between their
restriction to order $n$ and the fixed moment system $\langle X_n,\L
\rangle$. As part of the data, every moment system in $\MM_n$
includes a line bundle on the restriction $X_n$, so that we have a
natural forgetful functor $\MM_n \to \LL_n$.

\begin{lemma}\label{equi}
The forgetful functor $\MM_n \to \LL_n$ is an equivalence of
categories.
\end{lemma}

\proof{} Both $\MM_n$ and $\LL_n$ are stacks of groupoids on $X$ in
the Zariski topology, and the forgetful functor is compatible with
the stack structure. Therefore it suffices to prove the claim
locally on $X$. But for affine schemes $X$, the claim follows
immediately from Lemma~\ref{triv}.
\endproof

\proof[Proof of Theorem~\ref{main}.] Construct the desired formal
moment system as an inverse limit of order-$n$ moment systems, and
construct the order-$n$ moment systems by induction on $n$. The base
of the induction is the given Poisson scheme $X = X_0$. To establish
the induction step, note that If $H^2(X,\calo_X) = 0$, then every
line bundle on $X_{n-1}$ extends to a line bundle on
$X_n$. Moreover, if $H^1(X,\calo_X)=0$, then such an extension is
unique up to an isomorphism. Now apply Lemma~\ref{equi}.
\endproof

\proof[Proof of Proposition~\ref{conf}.] Consider the dual numbers
algebra $k\langle \eps \rangle = k[\eps]/\eps^2$, let $X' = X
\times_k k\langle\eps\rangle$, and consider two maps $p_0,p_1:X' \to
X$ -- $p_0$ is the natural projection, and $p_1$ is given locally by
$p_1^*(f) =f + \xi(f)\eps$. We can turn $X'$ into a Poisson scheme
in two ways: the obvious one, which we denote $X'_0$ (the map
$p_0:X' \to X$ is then Poisson), and the obvious one with the
Poisson bracket multiplied by $(1 + \lambda\eps)$, which we denote
$X'_1$ (the map $p_1:X'_1 \to X$ is Poisson).  Let $\X_i'$, $i =
0,1$ be the moment system associated by Theorem~\ref{main} to
$\langle X'_i,p_i^*L \rangle$.  Then since $H^1(X,\calo_X)=0$, the
line bundle $L$ admits an action of $\xi$, so that $p_0^*L\cong
p_1^*L$. We now notice that if we modify the moment system $\X'_0$
by composing the projection $\X_0' \to S$ with the map $S \to S$, $t
\mapsto \lambda t$, then it becomes a moment system for the pair
$\langle X_1',p_0^*L\cong p_1^*L\rangle$. By the unicity statement
of Theorem~\ref{main}, the two moment systems become isomorphic,
which proves the claim.
\endproof

\proof[Proof of Theorem~\ref{alg}.]
Assume that the Poisson scheme $X$ is algebraically convex, so that
we have the projective map $X \to Y = \Spec A$ with $A =
H^0(X,\calo_X)$. Take the canonical formal moment system $\X$
extending $X$ and its order-$n$ restrictions $X_n$. Let $A_n =
H^0(X_n,\calo_{X_n})$. Then since $H^1(X,\calo_X) = 0$, the
restriction maps $A_n \to A_{n-1}$ are surjective, and we have short
exact sequences
$$
\begin{CD}
0 @>>> t^nA @>>> A_n @>>> A_{n-1} @>>> 0.
\end{CD}
$$
Therefore $A_\infty = \displaystyle\lim_{\leftarrow}A_n$ is a
Noetherian algebra over $k$ which is flat over $k[[t]]$ and complete
with respect to the $t$-adic filtration. Since $H^2(X,\calo_X) = 0$,
an ample line bundle on $X$ extends to the formal scheme $\X$. Thus
it is a projective formal scheme over $\Y = \Spec A_\infty$, and we
can apply the Algebraization Theorem \cite[Th\'eor\`eme
  5.4.5]{EGA}. Moreover, by construction $A_\infty =
H^0(\X,\calo_X)$, so that $\X$ is indeed algebraically convex.

Consider now $\Tot(\L)$. The fact that the natural $\gm$-action is
Hamiltonian with the prescribed moment map is obviously precisely
equivalent to Definition~\ref{ms}~\thetag{ii}.  Assume that $X$ is
smooth and symplectic. By construction, we have smooth maps
$$
\Tot(\L) \to \X \to S
$$
whose fibers are, respectively, $\gm$ and $X$. Since $X$, $S$ and
$\gm$ are smooth, so is $\Tot(\L)$. It remains to prove that the
Poisson structure on $\Tot(\L)$ is non-degenerate at every point $z
\in \Tot(\L)$. The Hamiltonian vector field $H_t$ spans the subspace
in $T_z\Tot(\L)$ which is tangent to the $\gm$-orbit. Since $X$ is
symplectic, Hamiltonian vector fields $H_f$ span the tangent space
to $X$ at the image of the point $z$. These vector fields can be
lifted to $\Tot(\L)$ -- indeed, it suffices to lift the defining
function. All in all, we see that Hamiltonian vector fields span the
relative tangent subspace
$$
T_z(\Tot(\L)/S) \subset T_z\Tot(\L).
$$
But this is a subspace of codimension $1$, and it is easy to see
that the subspace spanned by Hamiltonian vector fields must be
even-dimensional. Thus they must span the whole $T_z\Tot(\L)$.
\endproof

\bigskip

\noindent
{\sc Steklov Math Institute\\
Moscow, USSR}\\
{\em E-mail\/}: {\tt kaledin@mccme.ru}


\begin{thebibliography}{EGA}

\bibitem[BGV]{BGV} N. Berline, E. Getzler, and M. Vergne, {\em Heat
kernels and Dirac operators}, Grundlehren der Mathematischen
Wissenschaften, {\bf 298}, Springer-Verlag, Berlin, 1992.

\bibitem[EGA]{EGA} A. Grothendieck, {\em \'El\'ements de
G\'eom\'etrie Alg\'ebrique, III}, Publ. Math. IHES {\bf 11}.

\bibitem[GK]{GK} V. Ginzburg and D. Kaledin, {\em Poisson
deformations of symplectic quotient singularities}, Adv. in Math,
{\bf 186} (2004), 1-57.

\bibitem[KV]{KV} D. Kaledin and M. Verbitsky, {\em Period map for
non-compact holomorphically symplectic manifolds}, GAFA {\bf 12}
(2002), 1265--1295.

\bibitem[K]{K} M. Kontsevich, {\em Deformation quantization of
algebraic varieties}, EuroConf\'erence Mosh\'e Flato 2000, Part III
(Dijon).  Lett. Math. Phys. {\bf 56} (2001), 271--294.

\end{thebibliography}
\end{document}